\documentclass{article}
\usepackage{amsfonts}
\usepackage[dvips]{graphics}
\usepackage{epsfig}
\usepackage{graphicx}
\usepackage{psfrag}
\usepackage{amsmath}
\usepackage{amssymb,subfigure}
\usepackage{natbib}
\usepackage{setspace}

\bibpunct{(}{)}{;}{a}{}{,}

\def\R{\mathbb R}

\def\epsilon{\varepsilon}
\def\ds{\displaystyle}

\newcommand{\be}{\par\nobreak\noindent  \begin{equation}}
\newcommand{\ee}{\end{equation}}
\newcommand{\baa}{\begin{array}}
\newcommand{\eaa}{\end{array}}
\newcommand{\ba}{\begin{eqnarray}}
\newcommand{\ea}{\end{eqnarray}}

\title{A population facing climate change: joint influences of Allee effects and environmental boundary geometry}

\author{Lionel Roques$^{\hbox{\small{ a,* }}}$, Alain Roques$^{\hbox{\small{ b }}}$, Henri Berestycki$^{\hbox{\small{ c }}}$, \\ and Andr\'e Kretzschmar$^{\hbox{\small{ a}}}$ \\
\\
\footnotesize{$^{\hbox{a }}$ Unit\'e Biostatistique et  Processus Spatiaux (BioSP), INRA,} \\
\footnotesize{ Domaine St Paul - Site Agroparc 84914 Avignon Cedex
9, France}\\
\footnotesize{$^{\hbox{b} }$ Station de Zoologie Foresti\`ere,
INRA,}\\
\footnotesize{ Av. de la Pomme de Pin, BP 20619, 45166 Olivet Cedex, France}  \\
\footnotesize{$^{\hbox{c} }$ Centre d'Analyse et de Math\'ematique
Sociales, EHESS,}\\
\footnotesize{54 bd Raspail 75270 Paris Cedex 06, France}  \\
 \footnotesize{$^{\hbox{*} }$ Corresponding author. E-mail: lionel.roques@avignon.inra.fr}
}

\date{}

\begin{document}

\maketitle

\begin{abstract}
As a result of climate change, many populations have to modify
their range to follow the suitable areas  - their ``climate
envelope" - often risking extinction. During this migration
process, they may face absolute boundaries to dispersal, because
of external environmental factors. Consequently, not only the
position, but also the shape of the climate envelope can be
modified. We use a reaction-diffusion model to analyse the effects
on population persistence of simultaneous changes in the climate
envelope position and shape. When the growth term is of logistic
type, we show that extinction and persistence are principally
conditioned by the species mobility and the speed of climate
change, but not by the shape of the climate envelope. However,
with a growth term taking an Allee effect into account, we find a
high sensitivity to the variations of the shape of the climate
envelope. In this case, the species which have a high mobility,
although they could more easily follow the migration of the
climate envelope, would be at risk of extinction when encountering
a local narrowing of the boundary geometry. This effect can be
attenuated by a progressive opening of the available space at the
exit of the narrowing,  even though it transiently leads to a
diminished area of the climate envelope.
\end{abstract}

{\bf Keywords:} Allee effect $\cdot$ Biodiversity $\cdot$ Climate
change $\cdot$ Climate envelope $\cdot$ Conservation $\cdot$
Mobility $\cdot$ Reaction-diffusion $\cdot$ Single species model

\section*{Introduction}\label{section1}

Over the last century, the global Earth's temperature has
increased by about 0.7$^\circ$C with, in the past 50 years, an
even faster warming trend, of about 0.13$^\circ$C per decade
\citep{IPCC2007}. The consequences of this warming on fauna and
flora are already visible, and well documented \citep[see
e.g.,][]{walt}. Especially, poleward and upward shifts of many
species ranges have been recorded, and are most likely linked to
the climate change \citep{py,pa}.

The  variations of  Earth's climate are in fact highly spatially
heterogeneous \citep{IPCC2007}. Thus, depending on the regions
they inhabit, species may be more or less subject  to high changes
in climate. Some species can adapt, whereas others, and especially
range-restricted species, risk extinction. A recent study based on
the species-area relationships, by \citet{thomas}, predicted
between $15\%$ and $37\%$ of species extinction in the next 50
years due to climate change, in sample regions covering $20\%$ of
the Earth's terrestrial surface.

For the next century, climate projections indeed predict a mean
increase of temperature of $1.8^\circ$C for minimum scenarios to
$4^\circ$C for maximum expected scenarios.

Recently, some authors proposed mathematical models for analysing
the factors that influence population persistence, when facing a
climate change. They used the notion of ``climate envelope'',
corresponding to the environmental conditions under which a
population can persist. They assumed that the conditions defining
this envelope did not change with time, while the envelope
location
 moved according to climate change
\citep{thomas}. Keeping $0.9^\circ$C/100~km as the poleward
temperature gradient, the above-mentioned temperature increases
should imply poleward translations of the climate envelope at
speeds of $2-4.5$~km/year.

A one-dimensional reaction-diffusion model has been proposed by
\citet{bdnz}. Results have been obtained, especially regarding the
links between population persistence and species mobility. Indeed,
if the individuals have a too low rate of movement, then the
population cannot follow the climate envelope and thus becomes
extinct. On the other hand, if the individuals have a very high
rate of movement, they disperse
 outside the climate envelope and thus the population also becomes
extinct. Other related works, based on one-dimensional
reaction-diffusion models, can be found in \citet{pl},
\citet{deasi}, \citet{pach}  and \citet{lutscher}.

The model that we consider here is derived from the classical
Fisher population dynamics model \citep[][]{fi,kpp}. In  a
two-dimensional bounded environment $\Omega$, the corresponding
equation is:
\par\nobreak\noindent  \begin{equation}\label{model}
\frac{\partial u}{\partial t}=D \nabla^2 u+u g(t,x,u), \ t\in
[0,+\infty), \ x\in\Omega\subset\R^2.
\end{equation}
The one-dimensional model considered by \citet{bdnz} is a
particular case of  (\ref{model}), with a logistic growth function
$u g(t,x,u)$. Indeed, the authors assumed the \textit{per capita}
growth rate $g$  to decrease with the population density $u$ for
all fixed $x$ and $t$.

Other types of growth  functions are of interest, especially those
taking account of Allee effect.  Allee effect occurs when the
\textit{per capita} growth rate reaches its peak at a strictly
positive population density. At low densities, the \textit{per
capita} growth rate may then become negative (strong Allee
effect). Allee effect is known in many species
\citep[see][]{allee,dennis,veit}. This results from several
processes which can co-occur \citep{berec}, such as diminished
chances of finding mates at low densities \citep{mc,oik}, fitness
decrease due to consanguinity  or decreased visitation rates by
pollinators for some plant species \citep{groom}. It is commonly
accepted that populations subject to Allee effect are more
extinction prone \citep{steph}.

In reaction-diffusion models, Allee effects are generally modelled
by equations of bistable type \citep[][]{fife,tur,shiallee}.
Mathematical analyses involving these equations have demonstrated
important effects of the domain's geometry, especially while
studying travelling wave solutions. These solutions generally
describe the invasion of a constant state, for instance where no
individuals are present, by another constant state (typically the
carrying capacity), at a constant speed, and with a constant
profile \citep[see][]{aw}. \citet{bh2} have proved that in an
infinite environment with hard obstacles, but otherwise
homogeneous, travelling waves solutions of bistable equations may
exist or not, depending on the shape of the obstacles. In the same
idea,  in an infinite homogeneous square cylinder,
\citet{chapuisat} proved that travelling wave solutions may not
exist if the cylinder's diameter is suddenly increased somewhere.
See also \citet{matano} for another related work.  In parallel,
\citet{keitt}, while studying invasion dynamics in spatially
discrete environments, showed that Allee effect can cause an
invasion to fail, and can therefore be a key-factor that
determines the limits of species ranges. In one-dimensional
models, Allee effect has also been shown to slow-down invasions
\citep{hurford}, or even to stop or reverse invasions in presence
of predators \citep{owen}, or pathogens \citep{hilker}. In a study
by \citet{tob}, empirical evidences have also been given regarding
the fact that geographical regions with higher Allee thresholds
are associated with lower speeds of invasion.

The aim of this work is to study, for the simple two-dimensional
reaction-diffusion model (\ref{model}), how  the population size
variation during a shift of the climate envelope depends on the
shift speed, the geometry of the environmental boundary, and on
the population mobility and growth characteristics.


The first section is dedicated to the precise mathematical
formulation of the model. Two types of growth functions are
considered, logistic-like, or taking account of an Allee effect,
with in both cases a dependence  with respect to the climate
envelope's position. We define three domain types, corresponding
to three kinds of geometry of the environmental boundaries. Domain
{\bf 1} is a straight rectangle, and domain {\bf 2} is  the union
of two rectangles of same width by a narrow corridor. The
comparison between domains {\bf 1} and {\bf 2} enables to analyse
the effects of a local narrowing of the habitat on population
persistence. As suggested by the work of \citet{chapuisat} and
\citet{bh2}, in the case with an Allee effect, extinction
phenomena may  not be caused directly by the reduction of climate
envelope due to the narrowing of the domain in the corridor, but
by its too sudden increase at the exit of the corridor. This is
why we introduced domains of type {\bf 3}, which correspond to the
union of two rectangles of same width by a narrow corridor, which
gradually opens over a trapezoidal region of length $h$ (the case
$h=0$ corresponds to domain {\bf 2}). Results of numerical
computations of the population size over $30$ years are presented
and analysed, under different hypotheses on the growth rate,
mobility, speed of climate change, and domain's shape. These
results are further discussed in the last section of this paper.

\section*{Formulation of the model}\label{section2}

The population dynamics is modelled by the following
reaction-diffusion equation:
\par\nobreak\noindent  \begin{equation*}
\frac{\partial u}{\partial t}=D \nabla^2 u+u g(t,x,u), \ t\in
[0,+\infty), \  x\in\Omega\subset\R^2.
\end{equation*}
Here, $u=u(t,x)$ corresponds to the population density at time $t$
and position $x=(x_1,x_2)$. The number $D>0$ measures the species
mobility, and $\nabla^2$ stands for the spatial dispersion
operator $\ds{\nabla^2 u=\frac{\partial^2 u}{\partial
x_1^2}+\frac{\partial^2 u}{\partial x_2^2}}$. The set $\Omega$ is
a bounded subdomain of $\R^2$. We assume reflecting  boundary
conditions (also called no-flux or Neumann boundary
conditions):\par\nobreak\noindent $$\frac{\partial u}{\partial
n}(t,x)=0, \hbox{ for }x\in \partial \Omega,$$where $\partial
\Omega$ is the domain's boundary and $n=n(x)$ corresponds to the
outward normal to this boundary. Thus, the boundary of the domain
$\Omega$, or equivalently the environmental boundary, constitutes
an absolute barrier that the individuals cannot cross.

\subsection*{Growth functions}\label{sub2.1}

The function $g$ corresponds to the \textit{per capita} growth
rate of the considered species. In our model, it can be of two
main types, $g_l$ or $g_a$. The first case  corresponds to a
logistic-like growth rate, depending on the position with respect
to the climate envelope: \be \left\{ \baa{ll}
g_l(t,x,u)=r^+(1-\frac{u}{K}), & \hbox{ if }x\in \mathcal{C}(t), \\
g_l(t,x,u)=r^+(1-\frac{u}{K})-r^-, & \hbox{ if }x \not \in
\mathcal{C}(t). \eaa \right. \label{gkpp} \ee The set
$\mathcal{C}(t)\subset \Omega$ corresponds to the climate envelope
at time $t$, the real number $r^+>0$ is the intrinsic growth rate
of the species inside the climate envelope, $K$ corresponds to the
carrying capacity inside the climate envelope, and $r^-$
corresponds to the drop in intrinsic growth rate, outside the
climate envelope.

In the second case,  a strong Allee effect is modelled. For fixed
values of $t$ and $x$, the function $g_a(t,x,u)$ does not attain
its maximum at $u=0$, and furthermore $g_a(t,x,0)<0$ for all
$t>0$, $x\in \Omega$. The typical form of \textit{per capita}
growth term taking account of an Allee effect is
\par\nobreak\noindent $$g_{a}=r\left(1-\frac{u}{k}\right)\left(\frac{u-c}{k}\right),$$where
$r$ is a growth term, $k$ is the environment's carrying capacity
and $c$ is the Allee threshold \citep[see e.g.,][]{lk,keitt}. For
the comparison with the logistic case (\ref{gkpp}) to stand, we
choose $r$ and $k$ such that $g_a(t,x,K)=g_l(t,x,K)=0$ for all
$x\in \mathcal{C}(t)$ and \par\nobreak\noindent $$\ds{\max_{u
\in(0,K), x \in \mathcal{C}(t)}g_a(t,x,u)=\max_{u \in(0,K), x \in
\mathcal{C}(t)}g_l(t,x,u)=r^+}.$$Moreover, we impose an ``Allee
threshold'' $c$ equal to $\rho K$, for $\rho\in(0,1)$, which means
that $g_a(t,x,u)<0$ as soon as $\ds{u<\rho K}$, and for all $x\in
\Omega$. Such conditions yields $k=K$ and
$r=\ds{\frac{4r^+}{(1-\rho)^2}}$. As in the logistic case, we
assume that $g_a$ drops by $r^-$ outside the climate envelope.
Finally, we obtain, \be \left\{ \baa{ll}
g_a(t,x,u)=r^+\frac{4}{(1-\rho)^2} \left(1-\frac{u}{K}\right)\left(\frac{u}{K}-\rho\right)& \hbox{ if }x\in \mathcal{C}(t), \\
g_a(t,x,u)=r^+\frac{4}{(1-\rho)^2}
\left(1-\frac{u}{K}\right)\left(\frac{u}{K}-\rho\right)-r^- &
\hbox{ if }x \not \in \mathcal{C}(t). \eaa \right. \label{gallee}
\ee

\begin{figure}
\centering
\subfigure[]{%
\includegraphics*[width=6cm]{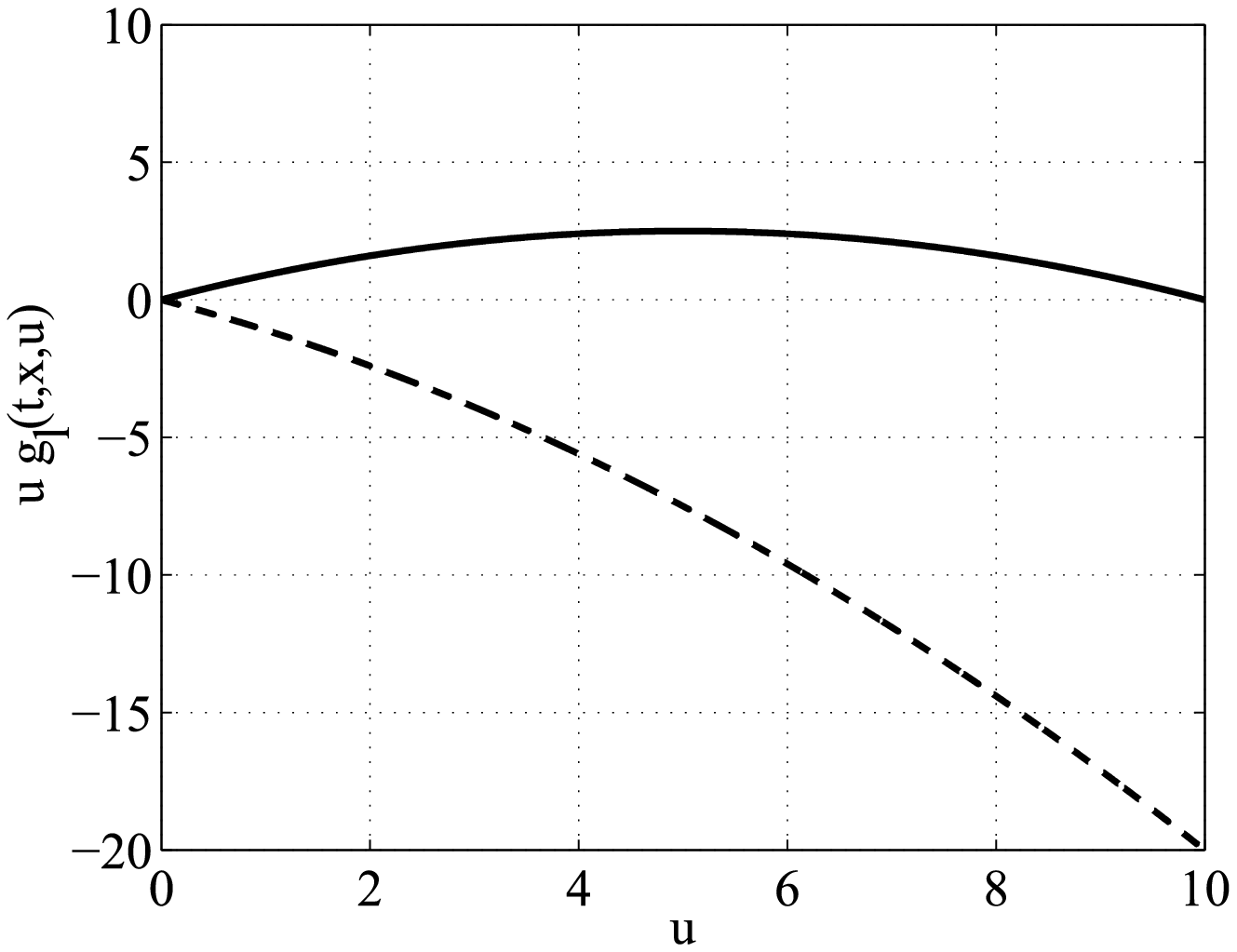}}
\subfigure[]{%
\includegraphics*[width=6cm]{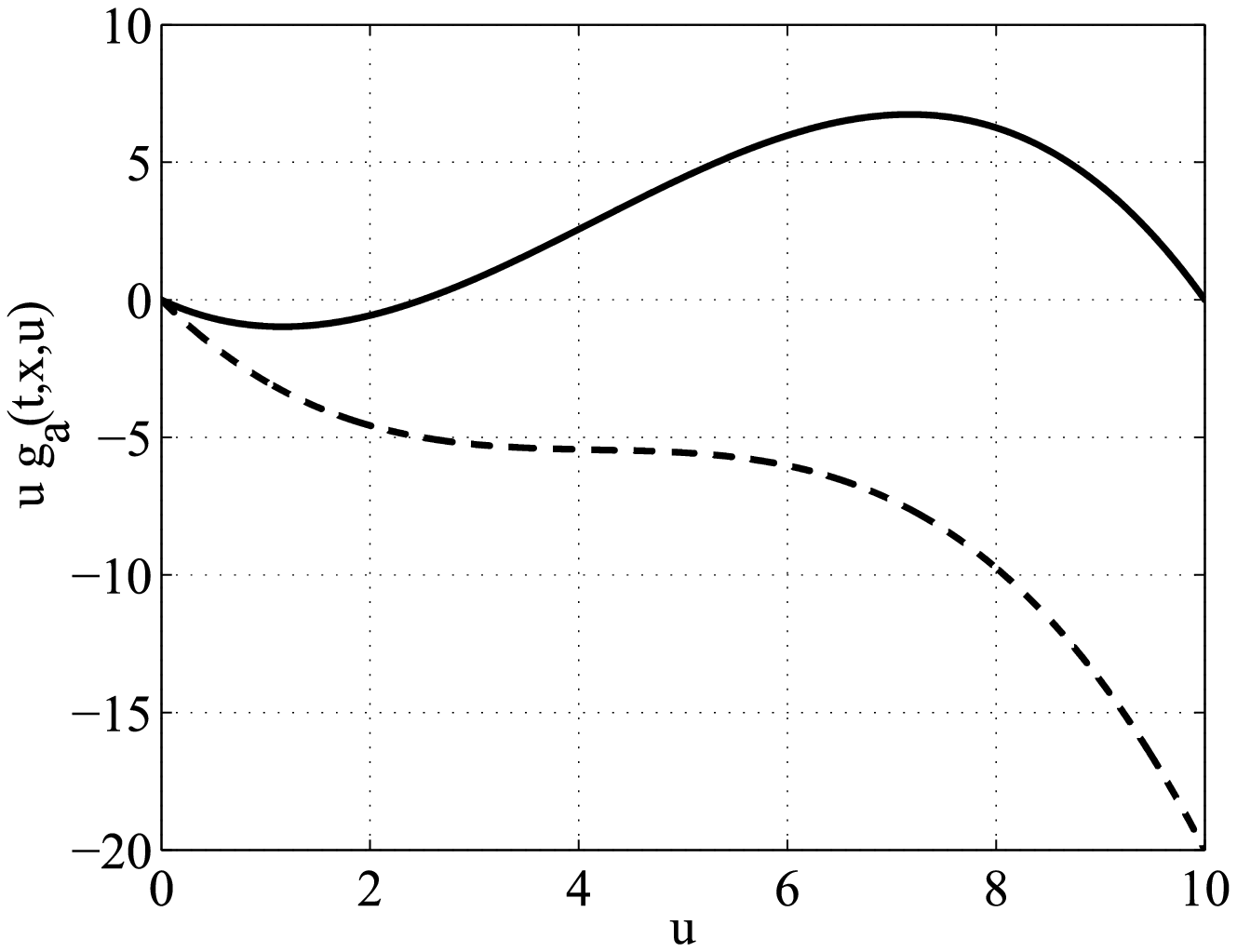}}
\caption{(a): profile of the growth rate $u\mapsto u g_l(t,x,u)$
inside the climate envelope ($x\in \mathcal{C}(t)$, \emph{solid
line}), and outside the climate envelope ($x \not \in
\mathcal{C}(t)$ \emph{dashed line}) in the logistic case. (b):
profile of the growth rate $u g_a(t,x,u)$ inside the climate
envelope (\emph{solid line}), and outside the climate envelope
(\emph{dashed line}), in the case with strong Allee effect. We
have set $r^+=1$, $r^-=2$, $K=10$, and
$\rho=0.25$.}%
\label{fig:gallee}
\end{figure}

The profiles of  the growth  functions $u g_l$ and $u g_a$ are
depicted in Figure \ref{fig:gallee}, for $x$ inside and outside
the climate envelope.

Throughout this paper, we make the hypothesis $r^+<r^-$.

\

\noindent {\bf Remark 1:} In both cases, the environment's
carrying capacity is equal to $K$ inside the climate envelope.
However, outside the envelope, the carrying capacity in the
logistic case would be $\ds{K\left(1-\frac{r^-}{r^+}\right)}$ for
$r^+ \geq r^-$. However, it is not defined when $r^+ < r^-$.
Indeed, in this last situation, the equation $u g_l(t,x,u)=0$ has
no positive solution. Similarly, in the case with Allee effect,
the equation $u g_a(t,x,u)=0$ has no positive solution  outside
the climate envelope as soon as $r^+ < r^-$, since $\ds{\max_{u
\in(0,K), x \in \mathcal{C}(t)}g_a(t,x,u)=\max_{u \in(0,K), x \in
\mathcal{C}(t)}g_l(t,x,u)=r^+}$. This means that the environment
is not suitable for persistence outside the climate envelope.

\

As explained in the introduction, we assume that the climate
envelope $ \mathcal{C}(t)\subset\Omega$ moves poleward, according
to the climate change. Assuming we are in the Northern hemisphere,
we consider this move to be of constant speed $v>0$, in the second
variable direction (to the ``North''). Furthermore, we assume $
\mathcal{C}(t)$ to be of constant thickness $L>0$:
\par\nobreak\noindent $$
 \mathcal{C}(t):=\{x=(x_{1},x_{2})\in \Omega, \hbox{ such that
 }x_{2}\in[vt,L+vt]\}.$$

\

\noindent {\bf Remark 2:} In this framework, the spreading speed
of the population, if it survives, is constrained by the speed at
which the climate envelope moves. Let us recall that, in an
homogeneous one-dimensional environment, with a logistic growth
rate $u r^+ (1-u/K)$, the spreading speed is $c^*_l=2\sqrt{r^+ D}$
\citep[see e.g.,][]{aw}. In such an homogeneous environment, but
with the growth rate, $4 u r^+/(1-\rho)^2
\left(1-u/K\right)\left(u/K-\rho\right)$, taking account of an
Allee effect, the spreading speed becomes $c^*_a=2\sqrt{r^+ D}
\sqrt{2} (1/2-\rho)/(1-\rho)$ \citep[see][]{lk}. We observe that
$c^*_a<c^*_l$, and that $c^*_a$ decreases with $\rho$. Thus, the
spreading speed is slowed down by the Allee effect. The case
$\rho=0$ corresponds to weak Allee effect: the \emph{per capita}
growth rate exhibits a maximum not at $u=0$, but remains positive
at low densities. In that case $c^*_a[\rho=0]=c^*_l/\sqrt{2}$. The
cases $\rho\ge 1/2$, corresponding to very strong Allee effects,
lead to population extinction \citep[see][]{lk}.

\

\subsection*{Geometry of the environmental boundary}\label{sub2.3}

We work in piecewise $C^1$ bounded domains $\Omega$ of $\R^2$ of
 three different types:

{\bf 1} A straight rectangle.

{\bf 2} The union of two rectangles of same width by a narrow
passage $\omega$.

{\bf 3} The union of two rectangles of same width by a narrow
passage $\omega$ followed by a trapezoidal region of height $h$.

\begin{figure}
\centering
\subfigure[Domain {\bf 1}]{%
\label{fs1}
\includegraphics*[width=4cm]{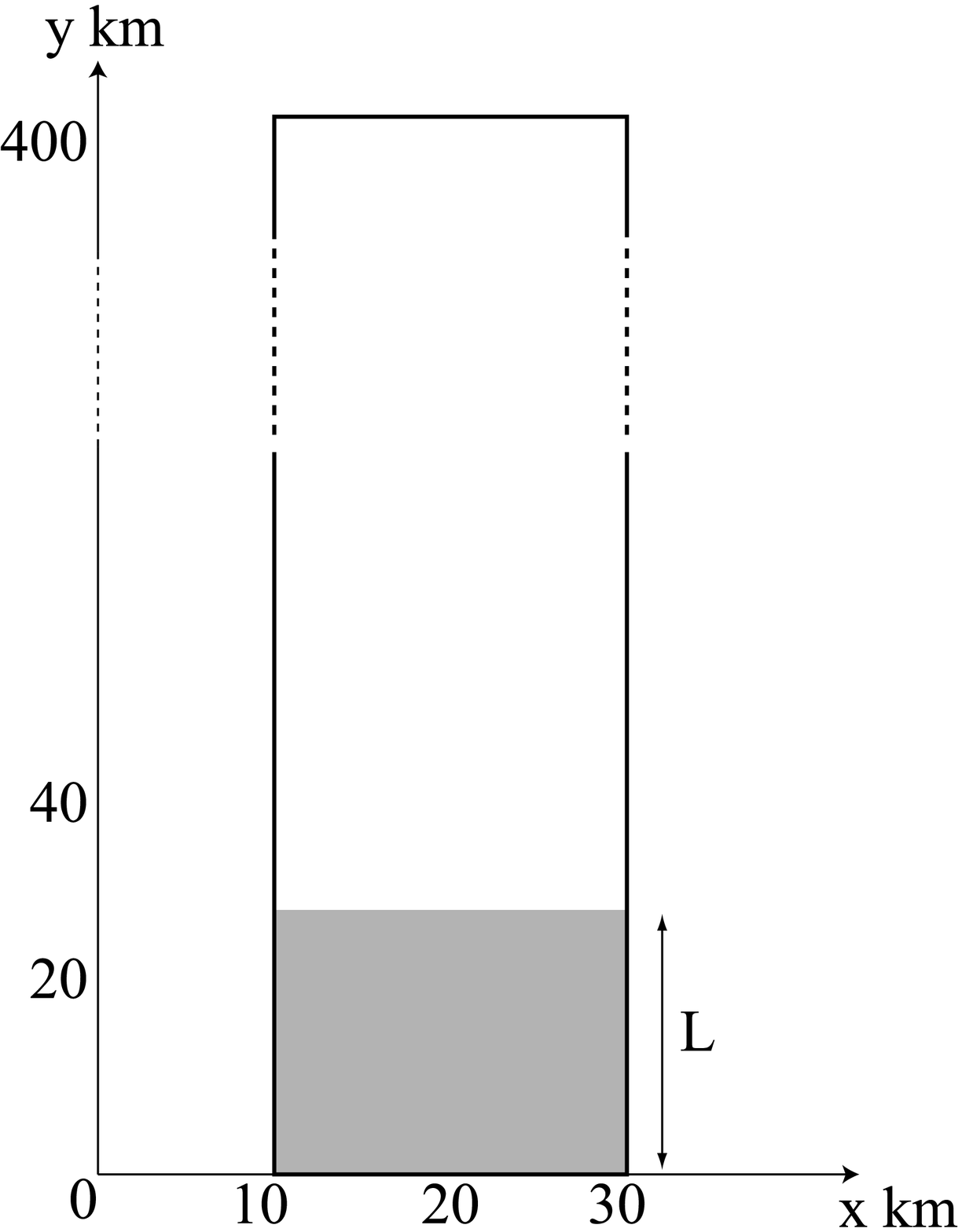}}
\subfigure[Domain {\bf 2}]{%
\label{fs2}
\includegraphics*[width=4cm]{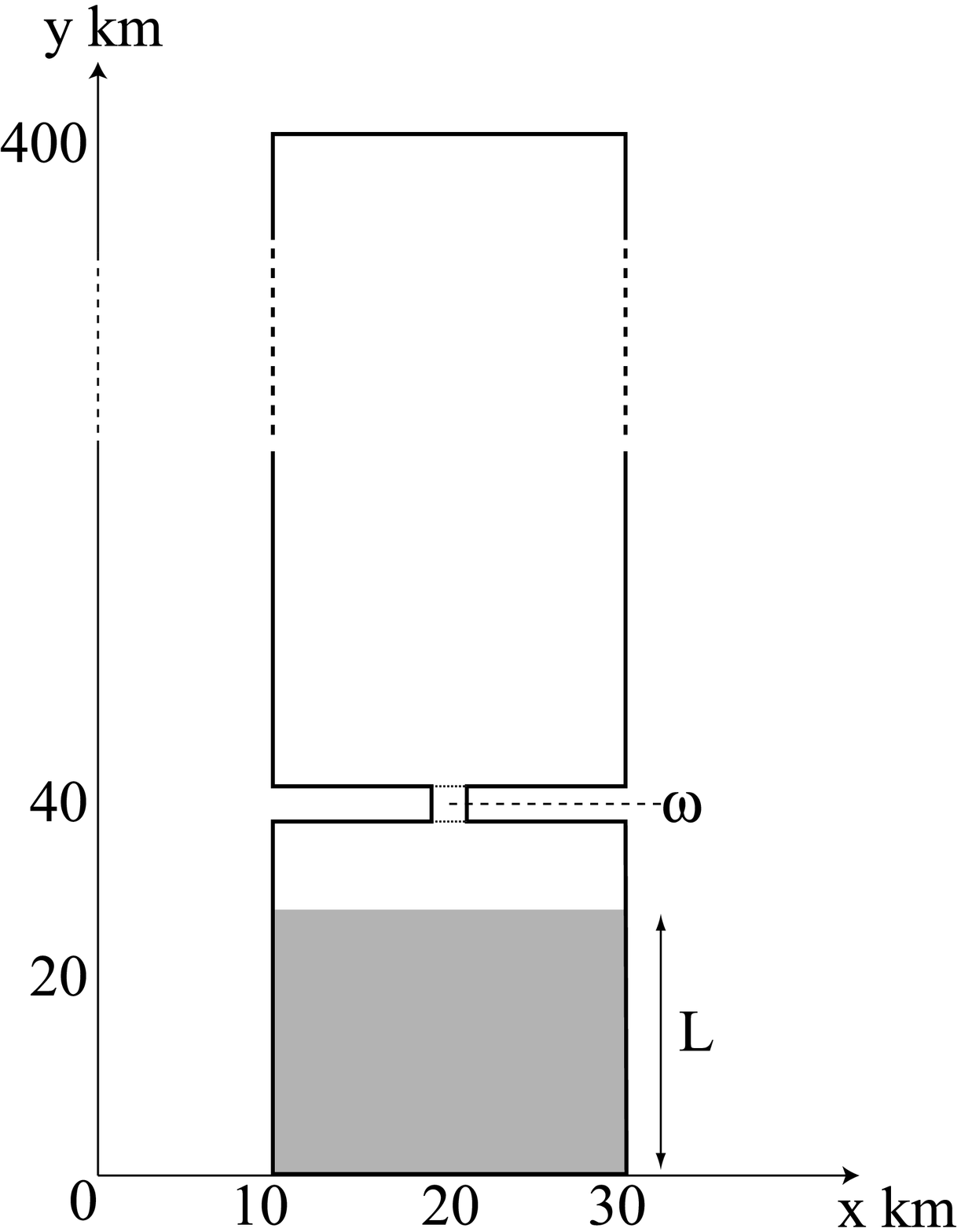}}
\\
\subfigure[Domain {\bf 3}]{%
\label{fs3}
\includegraphics*[width=4cm]{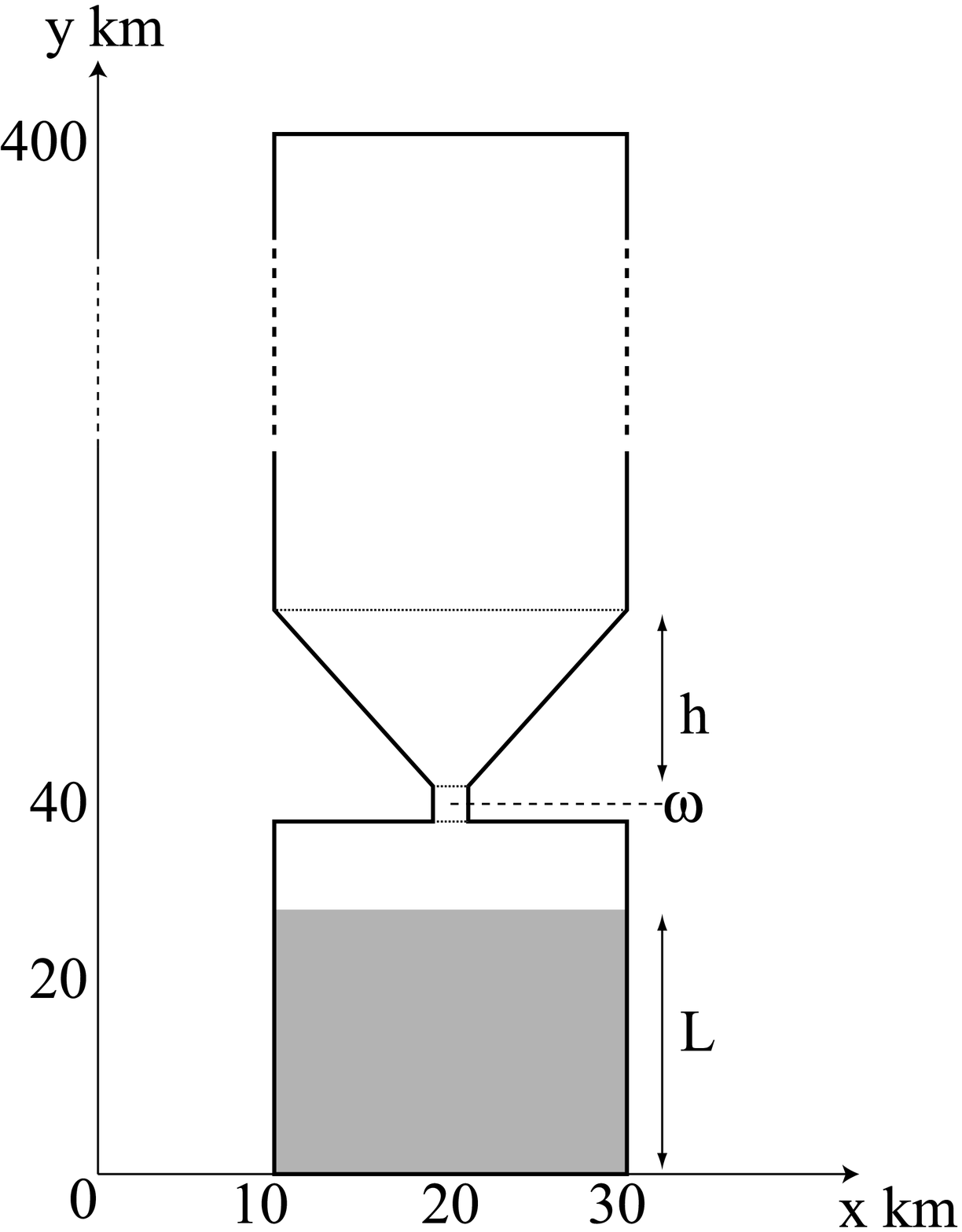}}
\caption{Three types of domains $\Omega$. The shaded area
corresponds to the climate envelope $\mathcal{C}(0)$ at $t=0$; it
is represented here with a length $L=30$km. The narrow passage
$\omega$ of domains {\bf 2} and {\bf 3}
 is separated by thin dotted lines, and is followed by a trapezoidal region of height $h$ in domain {\bf 3}.}%
\label{fig:dom}
\end{figure}

These three domain types, and their position relatively to the
origin, are depicted in Figure \ref{fig:dom}.

The domains are assumed to have a  width of $20$~km (on the larger
parts, for domains  {\bf 2} and {\bf 3}). In domains {\bf 2} and
{\bf 3},  the ``Southern" rectangles are assumed to have a length
of $40$~km. The width of the passage $\omega$ is assumed to be
$2$~km, and its length $4$~km.  The ``Northern" environmental
boundaries are situated far enough ($400$~km away from the
``Southern" boundaries), so that they have no influence on the
population. The climate envelope is assumed to have a latitudinal
range $L=30$~km.

\subsection*{Initial condition}\label{sub2.2}

We assume that the population was at equilibrium, and that the
climate envelope was stationary, before the considered period
starting at a time $t=0$. This means that the initial condition
$u(0,x)$ is a positive solution of the stationary equation:
\par\nobreak\noindent  \begin{equation}\label{eqsta}
-D \nabla^2 p=p \ g(0,x,p), \  x\in\Omega\subset\R^2,
\end{equation}
with $\frac{\partial p}{\partial n}(x)=0$ over $\partial \Omega$.
We place ourselves under the appropriate conditions for existence
of such positive stationary  solutions, in the logistic and Allee
effect cases (see the Appendix for more details).

\section*{Results}\label{section3}

Using a second order finite elements method, we computed the
solution $u(t,x)$ of the model (\ref{model}) with the growth
 functions and initial conditions discussed above. We focus here
on the population size, \par\nobreak\noindent $$
P(t):=\int_{\Omega}u(t,x)dx, $$ which was computed over $30$ years
in various situations.

Unless otherwise mentioned, we set $r^+=1$~year$^{-1}$ for the
intrinsic growth rate coefficient inside the climate envelope, and
we assumed that the \textit{per capita} growth rate was decreased
by $2$ outside the climate envelope: $r^-=2$~year$^{-1}$. We
assumed that $K=10$~individuals/km$^2$. In our computations, the
diffusion coefficient $D$ varied between $1$~km$^2$/year,
corresponding  to populations with low  mobility, e.g., some
insects species, and 50~km$^2$/year, corresponding to populations
with high mobility \citep[see][for some observed values of $D$,
for different species]{sk}. We first chose $\rho=0.25$, so that
the Allee threshold $\rho K$ is $2.5$~individuals/km$^2$. The
speeds used for the climate envelope shift varied between
$v=1$~km/year and $v=6$~km/year.

\subsection*{Population size over time}\label{sub_P(t)}

\begin{figure}
\centering
\subfigure[]{%
\includegraphics*[width=6cm]{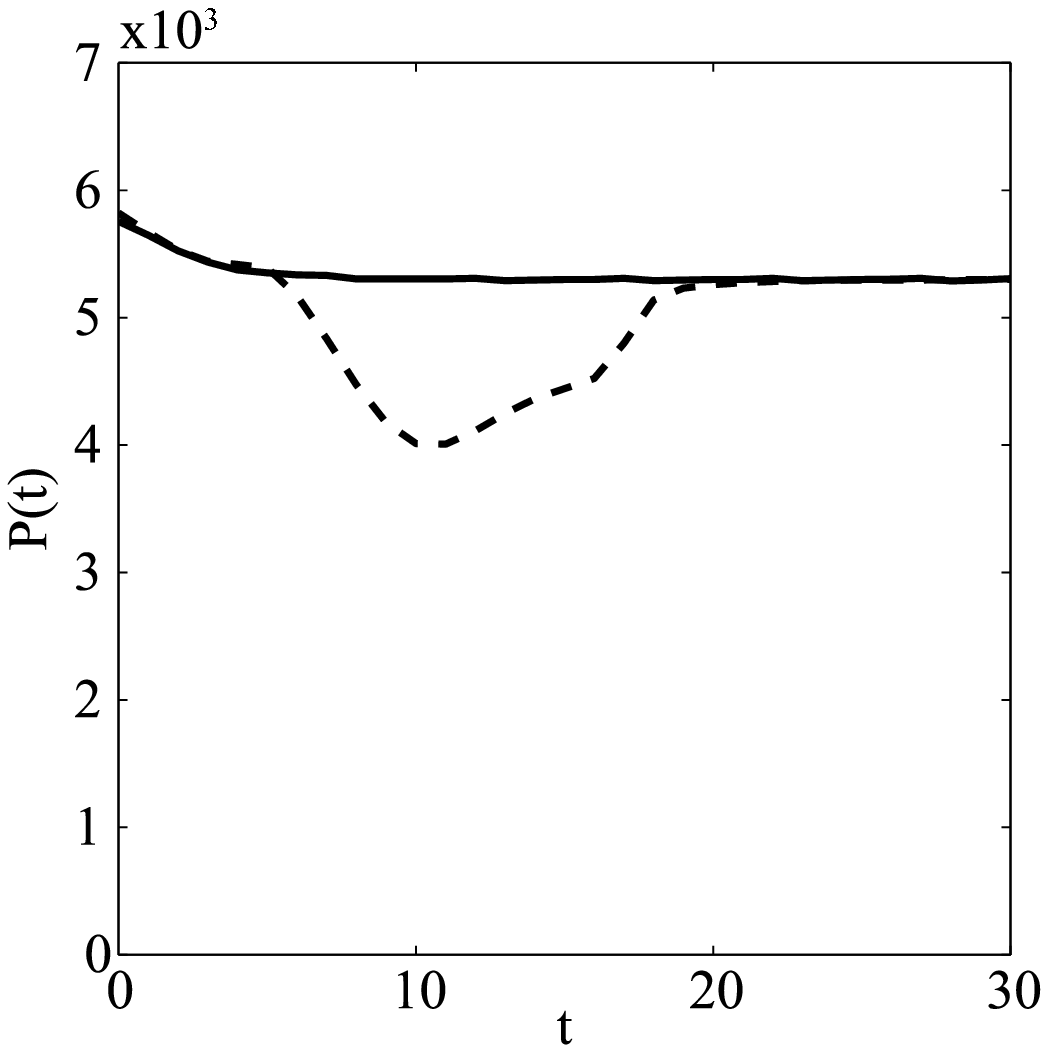}}
\subfigure[]{%
\includegraphics*[width=6cm]{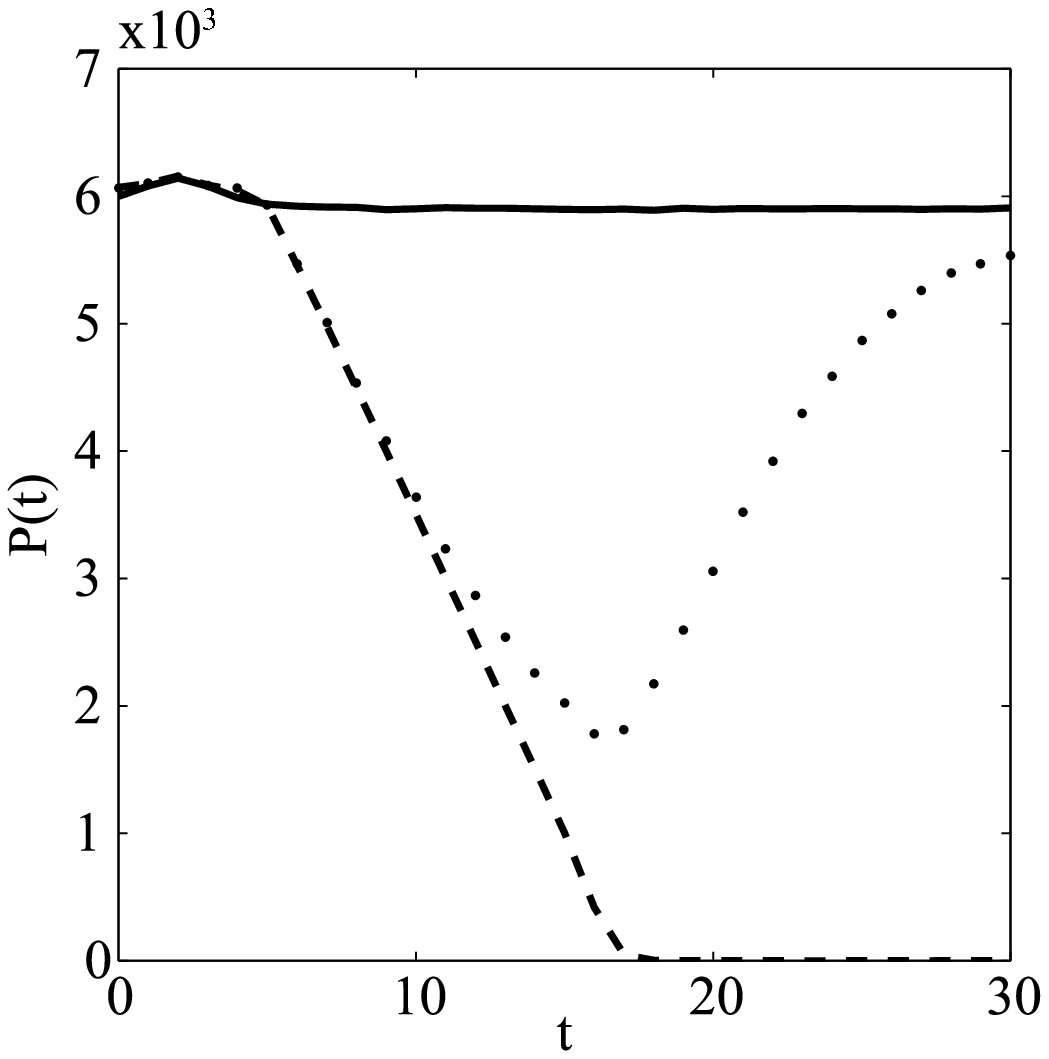}}
\caption{ Computation of the populations size $P(t)$ in function
of time $t$. (a) Case of the logistic growth rate $ug_l(t,x,u)$,
in domain {\bf 1} (\emph{solid line}), and in domain {\bf 2}
(\emph{dashed line}). (b) Case of the growth rate $ug_a(t,x,u)$
with Allee effect, in domains {\bf 1} (\emph{solid line}), {\bf 2}
(\emph{dashed line}) and {\bf 3} (\emph{dotted line}). For these
computation, we fixed $D=10$, $v=2.5$, $r^+=1$,
$r^-=2$, $K=10$, $\rho=0.25$ and $h=25$.}%
\label{fig:P(t)}
\end{figure}

In order to get a first  insight into the general behaviour of the
population size $P(t)$ for $t\in[0,30]$, we first computed it with
the fixed values $v=2.5$~km/year, and $D=10$~km$^2$/year. The
results are presented in Figure \ref{fig:P(t)}.

First notice that, as expected by our choice of  growth functions,
the initial population sizes, $P(0)$, are almost the same with
both growth  functions (\ref{gkpp}) and (\ref{gallee}): $5800$ and
$6000$ respectively in the logistic case and Allee effect case. A
slight transient increase or decrease of $P(t)$, probably due to
the initial effect of the Southern boundary, vanishes after
$5$~years.

In domain {\bf 1}, for both growth functions (\ref{gkpp}) and
(\ref{gallee}), after this short period of $5$~years, the
population sizes remain stable around their initial values.

On the contrary, the populations react very differently in the
domain {\bf 2}, depending on the type of their growth rate. In the
logistic case (Fig. \ref{fig:P(t)}a), the population size recovers
its initial value, after a transient decrease that lasts as long
as the corridor is included in the climate envelope. In the case
with Allee effect (Fig. \ref{fig:P(t)}b), the population size
declines to $0$ after $20$ years.

The behaviour of the population with Allee effect in domain {\bf
3}, with $h=25$~km (Fig. \ref{fig:P(t)}b), leads, as in domain
{\bf 2} with logistic growth, to the recovery of the initial level
of population.

As a preliminary conclusion of this first step, it comes that the
fate of the population will be driven by the interaction between
environmental parameters ($v$ and the type of domain) and
biological parameters ($D$, $r^+$, $r^-$, $K$, $\rho$,  and the
type of growth  function).

\subsection*{Intertwined effects of mobility and environmental parameters}\label{sub_Dv}

We computed the ``final" population size $P(30)$, for the range of
parameters $(v,D)\in (1,6)\times (1,50)$. We define the
``extinction region" as the portion of the parameter space leading
to $P(30)<1$.

\subsubsection*{Logistic case}\label{sub_Dv_log}

\begin{figure}
\centering
\subfigure[domain {\bf 1}]{%
\includegraphics*[width=5cm]{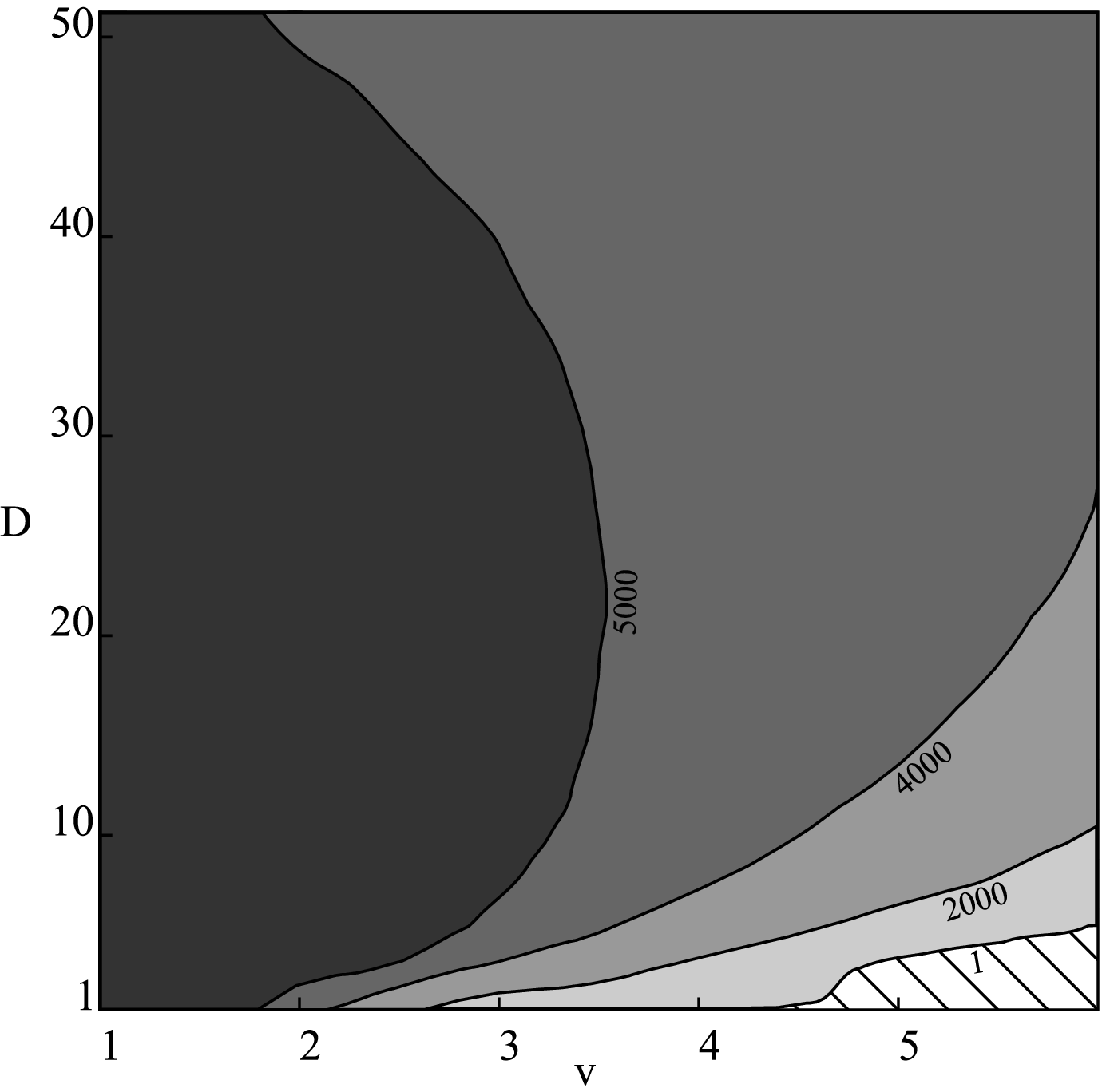}}
\subfigure[domain {\bf 2}]{%
\includegraphics*[width=5cm]{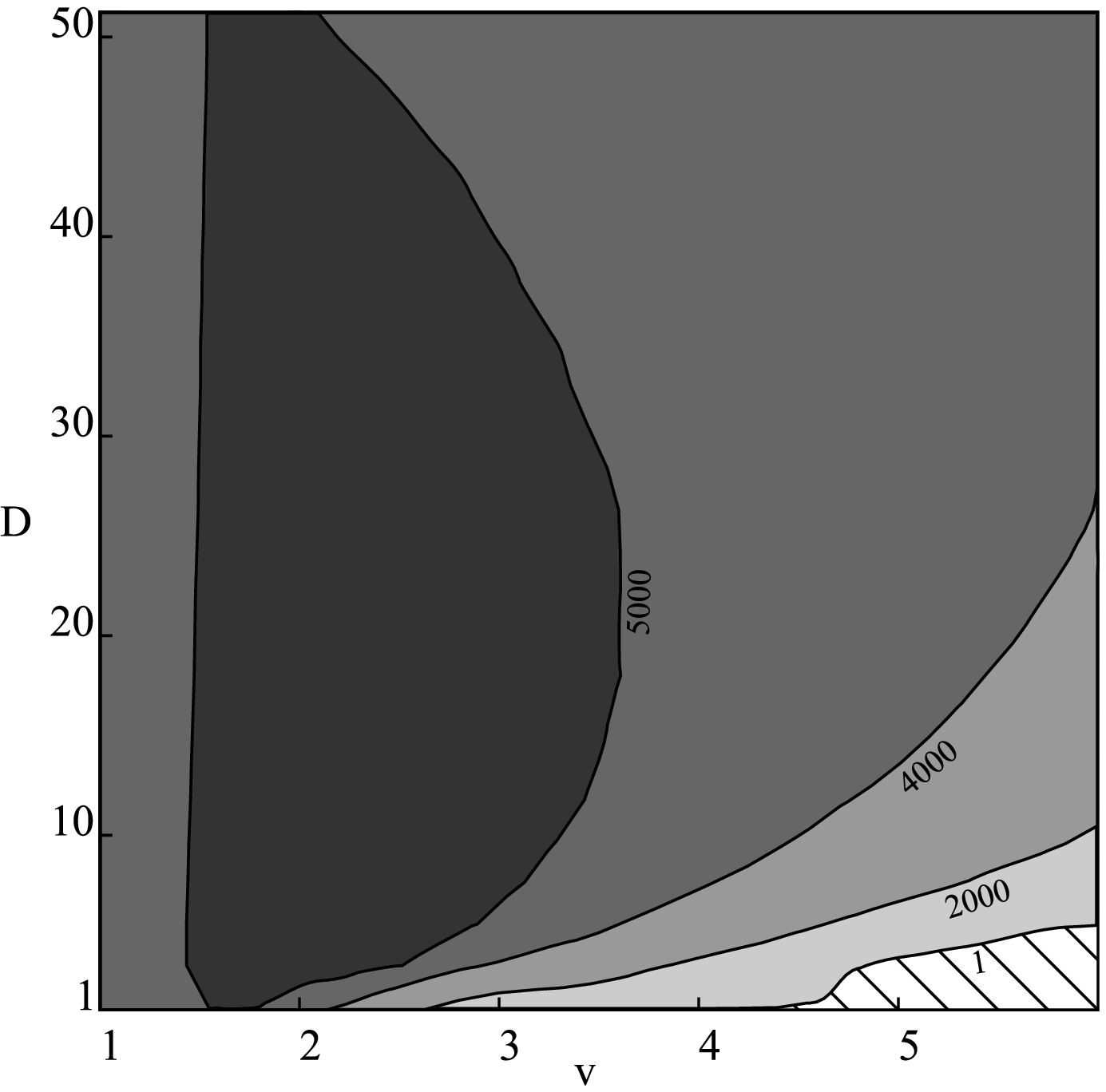}}
\caption{Logistic growth rate $ug_l(t,x,u)$. Contour lines of
$P(30)$ in function of the speed  of climate change $v$ and of the
diffusion coefficient $D$, in domains {\bf 1} (a) and {\bf 2} (b).
The hatched areas correspond to the
extinction regions. We fixed $r^+=1$, $r^-=2$ and $K=10$.}%
\label{fig:log2}
\end{figure}

The extinction region is very similar in domains {\bf 1} and {\bf
2} (Fig. \ref{fig:log2}), and corresponds to high values of $v$,
and low $D$ values. In those cases, the population mobility ($D$)
is not sufficient to follow the climate envelope (moving at speed
$v$). In domain {\bf 1}, for each value of $D$, $P(30)$ is
decreasing with respect to $v$, while in domain {\bf 2}, $P(30)$
is decreasing in $v$ for $v\ge 1.5$ (whenever $v< 1.5$, the
climate envelope has not totally crossed the corridor $\omega$ at
$t=30$). As already proved in the one-dimensional case
\citet{bdnz}, and for the same reasons recalled in the
introduction section, we here observe that population sizes
$P(30)$ are not monotonically linked with the parameter $D$.

\subsubsection*{Case with Allee effect}\label{sub_Dv_allee}

\begin{figure}
\centering
\subfigure[domain {\bf 1}]{%
\includegraphics*[width=5cm]{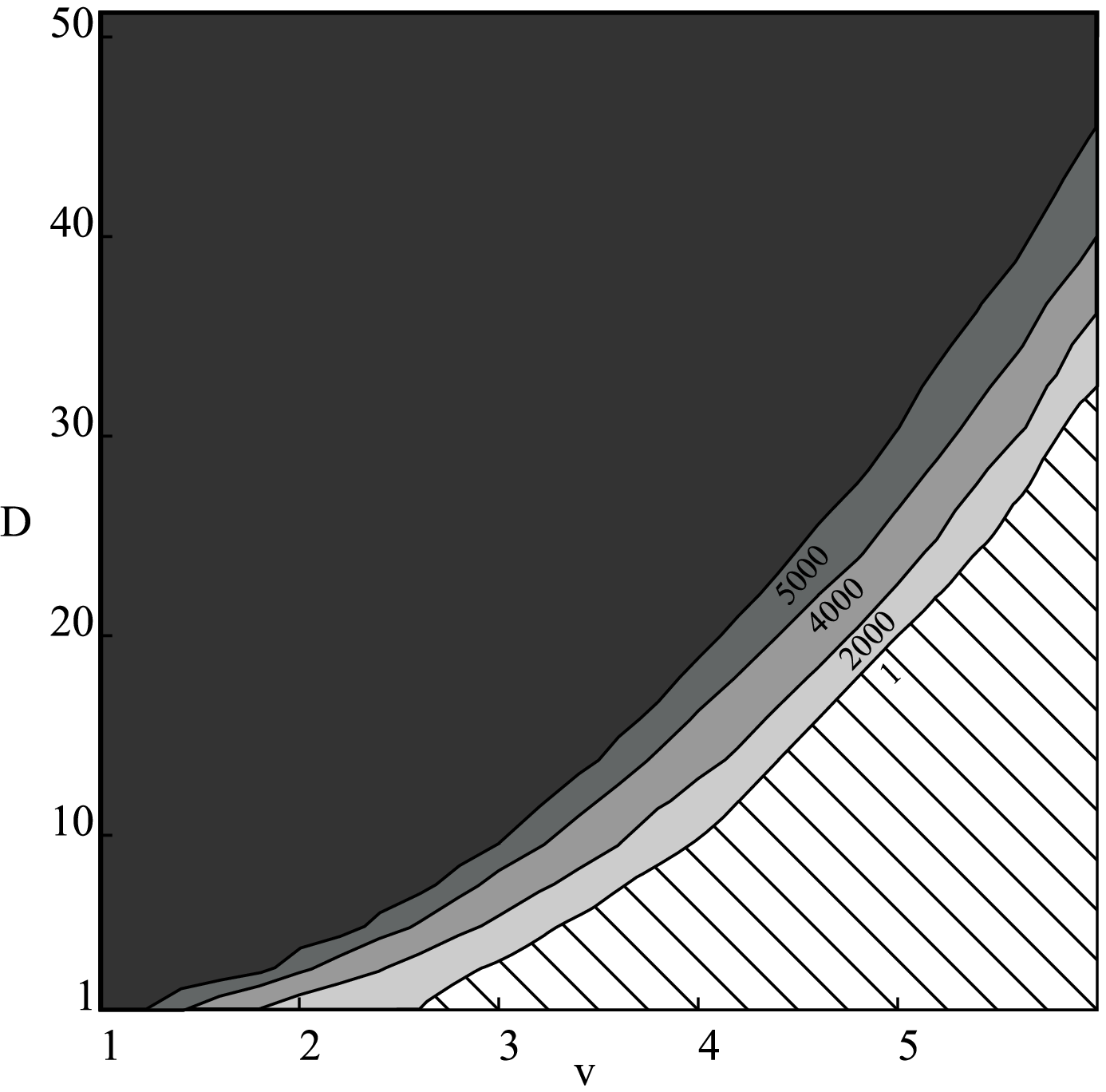}}
\centering
\subfigure[domain {\bf 2}]{%
\includegraphics*[width=5cm]{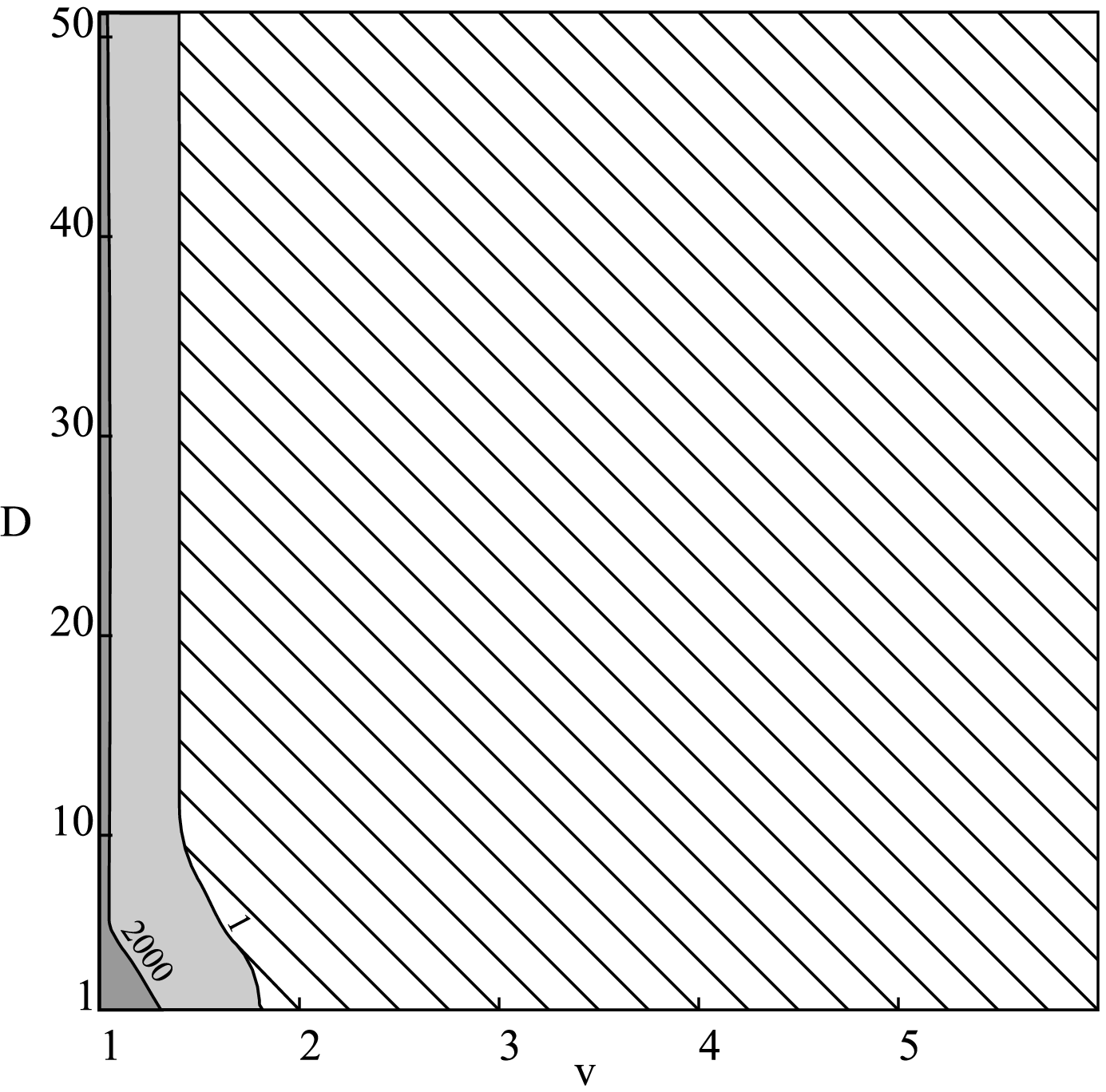}}
\subfigure[domain {\bf 3}]{%
\includegraphics*[width=5cm]{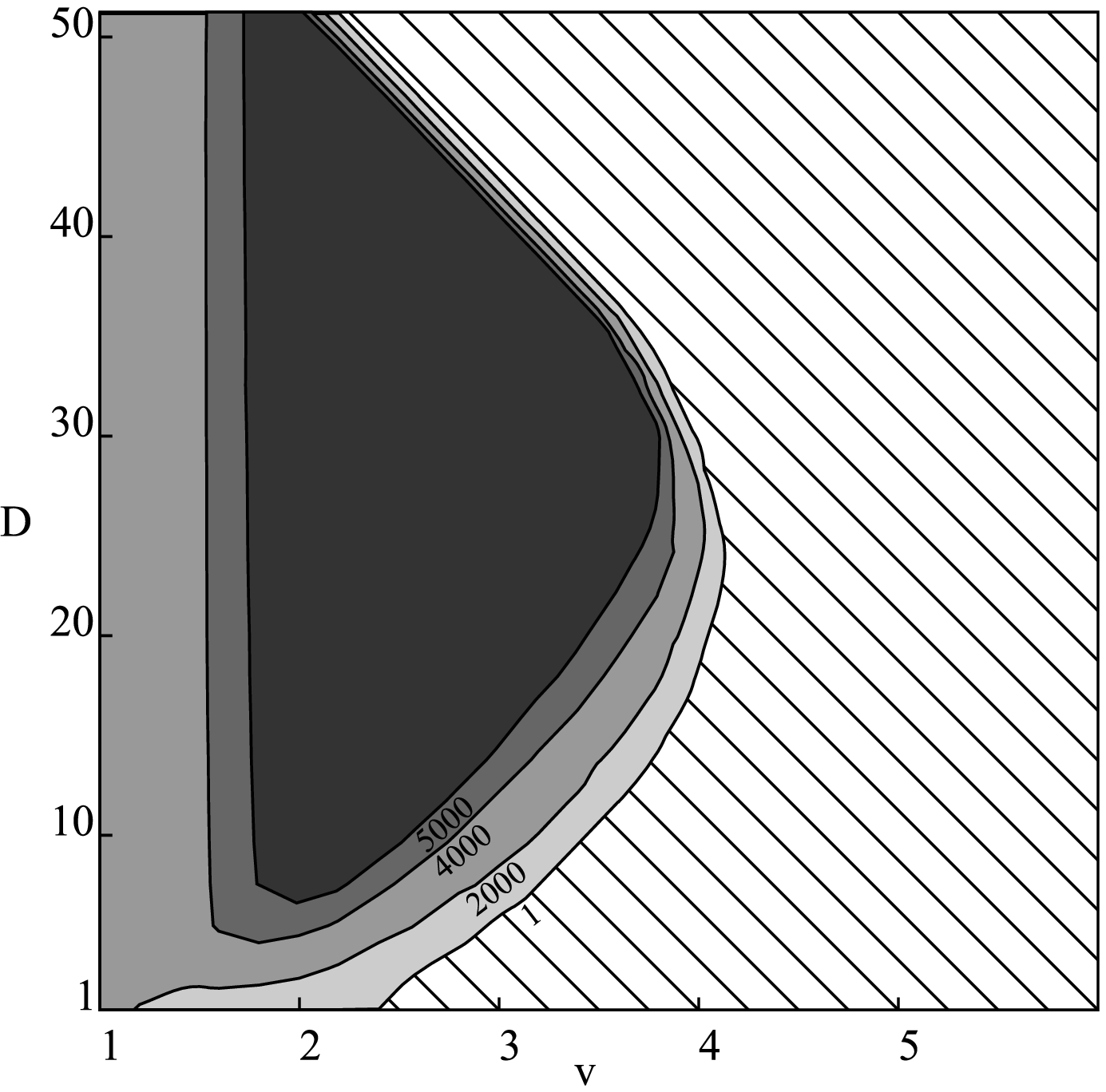}}
\caption{Case with Allee effect. (a), (b) and (c): contour lines
of $P(30)$ in function of $v$ and $D$ in domains {\bf 1}, {\bf 2}
and {\bf 3}, respectively. The hatched areas correspond to the
extinction regions. For these computations, we fixed $r^+=1$,
$r^-=2$, $K=10$, $\rho=0.25$ and
$h=25$.}%
\label{fig:allee2}
\end{figure}

In   domain {\bf 1} (Fig. \ref{fig:allee2}a), $P(30)$ is
decreasing with respect to $v$ and increasing with $D$. The
extinction region again corresponds to  high values of $v$, and
low $D$ values. However, the extinction region is wider than in
the logistic case. In domain {\bf 2} (Fig. \ref{fig:allee2}b), the
population does not survive as soon as $v>1.8$. The extinction
region is therefore very wide compared to the one of domain {\bf
1}.

Contrarily to the logistic case, a marked difference between the
fate of populations in domains {\bf 1} and {\bf 2} does exist.
This reflects the sensitivity to environmental boundary geometry,
of populations subject to an Allee effect. In domain {\bf 3}, with
$h=25$, where the opening at the exit of the corridor is more
progressive than in domain {\bf 2}, we observe that the extinction
region is reduced (Fig. \ref{fig:allee2}c). It remains however
larger than in domain {\bf 1}, especially for large values of $D$.

\subsection*{Allee effect and boundary geometry}\label{sub_Dh}

\subsubsection*{Effect of gradual boundary change at the exit of the corridor}\label{sub_Dh_1}

\begin{figure}
\centering
\includegraphics*[width=6cm]{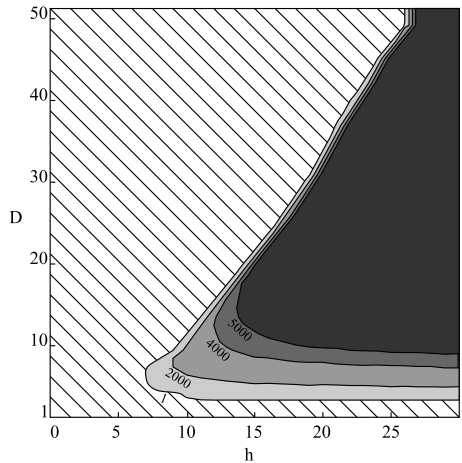}
\caption{Case with Allee effect. Contour lines of $P(30)$ in
domain {\bf 3}, in function of $D$ and of the height $h$ of the
trapezoidal region. The hatched area corresponds to the extinction
region. The other parameters values were
$r^+=1$, $r^-=2$, $K=10$, $\rho=0.25$ and $v=2.5$.}%
\label{fig:cas3h}
\end{figure}

We have seen, in the case of population subject to an Allee
effect, that domain {\bf 3}, with $h=25$, gives population
persistence more chances than domain {\bf 2} (which corresponds to
domain {\bf 3}, with $h=0$). Let us fix $v=2.5$, and see how
persistence depends on $h$ and $D$. The results of the computation
of $P(30)$,  for $(h,D)\in [0,30]\times(1,50)$ are presented in
Figure \ref{fig:cas3h}.

For each value of $D$, $P(30)$ increases with $h$; thus higher $h$
values give population persistence more chances. Furthermore, for
each $D\ge 3$, we observe that we can define a real number $h^*$
such that the population goes to extinction if $h < h^*$,  and the
population survives if $h > h^*$. Remarkably,  $h^*$ defines two
types of geometry, characterised by the progressiveness of the
opening at the exit of the corridor $\omega$, which lead to
extinction or survival, respectively. For most values of $D$
($D\ge 5$), $h^*$ linearly increases with $D$. For small values of
$D$ ($D<3$), extinction occurs independently of the value  $h$,
which simply means that, as in domain {\bf 1}, the population
cannot follow the climate envelope shift.

These results indicate that the  increased extinction risk in
domain {\bf 2}, compared to domain {\bf 1}, may not be caused by
the shrinking of the climate envelope and its fragmentation in two
parts, but by the lack of progressiveness of the opening at the
exit of the corridor.

\subsubsection*{Population density at the exit of the corridor governs the fate of the population}\label{sub_Dh_2}

\begin{figure}
\centering
\includegraphics*[width=6cm]{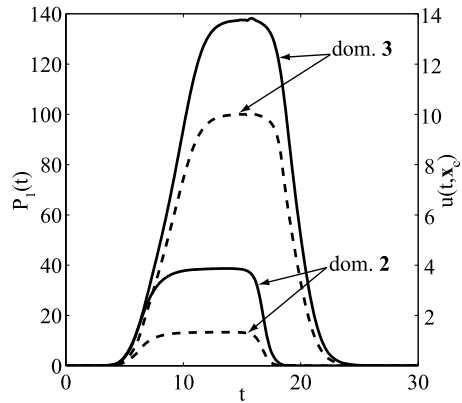}
\caption{ Case with Allee effect. \emph{Solid lines}: Number of
individuals $P_1(t)$ in the region $\Omega_1=\Omega\cap \{y\in
(44,48)\}$, corresponding to the four kilometres following the
corridor, in domains {\bf 2} and {\bf 3}. \emph{Dashed lines}:
population density $u(t,c)$ at the ``centre"
$\mathbf{x}_c=(10;46)$ of $\Omega_1$. For this computation, we
fixed $D=10$, $r^+=1$, $r^-=2$, $K=10$, $\rho=0.25$ and $v=2.5$.}%
\label{fig:omega1}
\end{figure}

Let us focus on what happens at the exit of the corridor $\omega$.
We considered the region $\Omega_1=\Omega\cap \{y\in (44,48)\}$,
corresponding to the four kilometres following the corridor.
 We computed the number $P_1(t)=\ds{\int_{\Omega_1}u(t,x)dx}$ of
individuals in this region, and the population density $u(t,c)$ at
the ``centre" $\mathbf{x}_c$ of $\Omega_1$, of coordinates
$(10;46)$.
 We compared domains {\bf 2} and {\bf 3}, with
the fixed parameters values $v=2.5$, $D=10$, $\rho=0.25$, and
$h=25$ (which is larger than $h^*=10$ in this case, see Fig.
\ref{fig:cas3h}). The results are presented in Figure
\ref{fig:omega1}.

We observe that, at the beginning of the invasion of $\Omega_1$,
for $t\in[4,7]$, the populations $P_1(t)$ are almost identical in
the two domains. Thus, the boundary geometry after the corridor
has no effect on the number of individuals which first invade the
region $\Omega_1$. On the contrary, during the same period of
time, the population density at $\mathbf{x}_c$, which is situated
at the exit of the corridor, remains lower in domain {\bf 2}, and
never reaches the Allee threshold $\rho K=2.5$. Thus, in domain
{\bf 2}, the dispersion of the few invaders into an open wide
space, after the corridor, leads to a low population density, and
therefore to a negative value of the \textit{per capita} growth
rate $g_a(t,x,u)$ (it would be false for a logistic growth
function).  The higher density in domain {\bf 3}, for $h> h^*$,
results from the progressiveness of the domain width offered to
population development.  Even if the individuals are scarce in
$\Omega_{1}$, the population density at  $\mathbf{x}_c$ exceeds
the Allee threshold at $t=7$.
 That makes the population
increase, leading to survival (survival can indeed be observed in
Fig. \ref{fig:allee2}c). Thus, the fate of the population seems to
depend on whether or not the population density at the exit of the
corridor reaches the Allee threshold.

\begin{figure}
\centering
\includegraphics*[width=6cm]{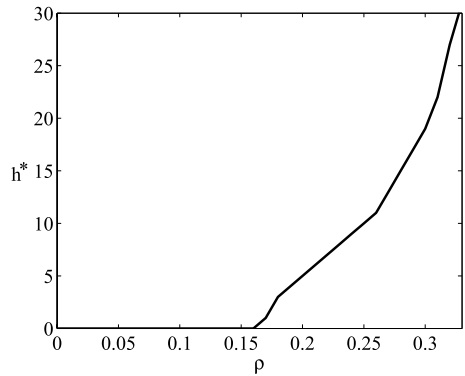}
\caption{Case with Allee effect. Minimal value $h^*$ of $h$
required for persistence, in terms of the parameter $\rho$. The
other parameters values were $D=10$, $r^+=1$, $r^-=2$,
$K=10$ and $v=2.5$.}%
\label{fig:rho}
\end{figure}

This leads us to analyse the effect of the Allee threshold. For
the same values $v=2.5$ and $D=10$, the function $\rho \mapsto
h^*$ is presented in Fig. \ref{fig:rho}. For small values of
$\rho$ ($\rho\le 0.16$), we observe that  $h^*=0$, i.e., the
population even survives in domain {\bf 2}. Then, $h^*$ increases
with $\rho$. For higher values of $\rho$, the population never
persists, at least  when $h\leq 30$ (see also the end of
remark~2). Finally, the stronger the Allee effect, the more
progressive must be the opening at the exit of the corridor
$\omega$.

On the contrary, the degree of hostility outside $\mathcal{C}(t)$,
quantified by $r^-$, does not seem to affect much $h^*$. Indeed,
computing $h^*$, for  $v=2.5$, $D=10$ and $\rho=0.25$, we always
obtained $h^*=10$, for $r^-\in (1,2)$ (remember that, in this
work, it is assumed that $r^->r^+$, see remark~1).

\section*{Discussion}\label{section4}

Using a  two-dimensional reaction-diffusion model, we have studied
the fate of populations with different mobility and growth
characteristics, facing environmental changes: a shift of their
climate envelope, the shape of which can be modified by
environmental boundary geometry.

The growth  functions we considered were of two types, logistic,
with a \textit{per capita} growth rate modelled by (\ref{gkpp}) or
taking account of a strong Allee effect (\ref{gallee}). Boundary
geometry diversity has been summarised into three schematic domain
shapes (Fig. \ref{fig:dom}).

In the logistic case, the other biological parameters being fixed,
the response (survival or extinction) to a climate envelope shift
at speed $v$ is simply determined by the mobility of the
population, which is measured by the diffusion parameter $D$. The
minimum required mobility for survival increases with the speed
$v$. This response is independent of the environmental boundary
geometry. This clearly appeared comparing domains {\bf 1} and {\bf
2}. In fact, we conjecture that it is true for any local
perturbation of the boundary (it could be proved using the methods
of Berestycki and Rossi, in press).

In the case with Allee effect, we observed a more complex pattern
of interactions among biological and environmental parameters. In
the straight domain {\bf 1}, similarly to the logistic case, the
higher the climate envelope shift speed, the higher the required
mobility for survival. Additionally, this required mobility, at
any $v$, is always higher than  in the logistic case.

Facing a local narrowing of the environment, as in domain {\bf 2},
contrarily to the logistic case,  the chances of survival
dramatically decrease in the presence of an Allee effect. In that
case, high mobility is not sufficient to face  climate envelope
shifting at high speeds. Indeed, the population density drop at
the exit of the corridor exhibits the sensitivity of this type of
populations to low densities. This drop can be attenuated by a
progressive opening of the available space at the exit of the
corridor, as in domain {\bf 3}. Even though this geometry
transiently leads to a diminished area of the climate envelope, it
finally results in higher chances of survival. The opening has to
be all the more progressive since the Allee threshold and the
species mobility are high.

Thus, a population subject to an Allee effect should have a
mobility which is neither too low, in order to be able to follow
the climate envelope, nor too high, in order to overcome various
changes of the shape of the climate envelope. For this reason,
species having an increasing dependence of their mobility with
respect to their density should be more robust to climate change
\citep[see][for a survey of reaction-diffusion models with
density-dependent dispersal terms]{oku}. On the other hand, under
the environmental changes considered in this paper, populations
with negatively density-dependent dispersal, such as in the  model
proposed by \cite{kmcc} for the dispersal of early Palaeoindian
people in North America, should have high probabilities of
extinction, if an Allee effect was assumed.

The reflective boundary conditions that we used throughout this
paper mean that individuals encountering the environmental
boundary are reflected inside the domain. These boundary
conditions can be encountered in many real-world situations,
corresponding to cliffs, rivers or coasts \citep{sk,fahrig2004}.
Other boundary conditions could have been considered. With Robin
boundary conditions for instance, a part of the individuals
crosses the boundary. These boundary conditions write $\partial
u/\partial n+\varepsilon u=0$ \citep[see e.g.,][for some details
on this type of boundary conditions]{ccL}. Preliminary numerical
computations have shown that our results still hold in such a
case, at least for small positive values of $\varepsilon$,
corresponding to few individuals crossing the environmental
boundary.

It was stated in \citet{pa}, on the basis of empirical studies,
that range-restricted species, like polar and mountain-top
species, were at high risk of extinction induced by warming. Our
study suggests that some other species may fail to expand poleward
and that their capacity to expand is linked to the geometry of the
geographical limits. Apart from laboratory tests under controlled
conditions, the diversity of arrangements of the Alpine valleys
could allow us to see whether the results of this paper can be
observed in natural conditions. Indeed, insects such as pine
processionary moth, the present range expansion of which is
undoubtedly related to climate change \citep{rob}, may allow us to
test statistically the relationships between insect progression
and the geometry of the Alpine corridors.

In this paper, a southern retraction of the climate envelope was
assumed. This can be directly linked, for some species, to the
fact that they are sensitive to high temperatures. This can also
be caused indirectly by competition with other species, the  range
of which is shifting to the North. For some species, however, this
retraction does not occur, leading to an expansion of the species
range. This could be easily integrated in our model, by setting
\par\nobreak\noindent $$
 \mathcal{C}(t):=\{x=(x_{1},x_{2})\in \Omega, \hbox{ such that }x_{2}\in[0,L+vt]\}.
$$ In this case, comparable results should be obtained, with
stagnation instead of extinction.

\section*{\protect Appendix: initial conditions}

When the growth rate is of logistic type, with $g=g_l$ satisfying
(\ref{gkpp}), a sufficient condition for the existence of a
positive solution $p_{l}$ of (\ref{eqsta}) can be derived by
finding an appropriate sub-solution of (\ref{eqsta}) \citep[see
e.g.,][]{am2}; indeed, the constant $K$ is readily a
super-solution. With our choices of $\Omega$, $\mathcal{C}(0)$ and
boundary conditions, for $\alpha>0$ small enough, the function
\par\nobreak\noindent $$\phi(x,y)=\alpha \sin\left(\frac{\pi}{2}+\frac{\pi y}{2
L}\right)\hbox{ for }y<L, \hbox{ and }\phi(x,y)=0 \hbox{ for
}y\geq L,$$is such a sub-solution, as soon as $\frac{D \pi^2}{4
L^2}<r^+.$ Under this condition, the function  $p_{l}$ is then the
unique positive and bounded solution of (\ref{eqsta}) \citep[it
can be proved as in][]{bhr1}. In the logistic case, we assume that
$u(0,x)=p_{l}$.

In the case (\ref{gallee}) with Allee effect, the condition for
the existence of a solution of (\ref{eqsta}), with $g=g_a$ is more
complex, and multiple solutions may exist. However, the existence
of a positive and bounded solution of (\ref{eqsta}) is still
granted when $\mathcal{C}(0)$ contains a sufficiently large ball
$B_{R}$. Indeed, for $R$ large enough, there exists a positive
solution $\psi$ of (\ref{eqsta})  on $B_{R}$, with Dirichlet
boundary conditions  \citep{bl}. The function $\psi$, extended by
$0$ to $\Omega$, is then a sub-solution of (\ref{eqsta}); the
constant $K$ is again a super-solution; this implies the existence
of a solution $p_a$ of (\ref{eqsta}) on $\Omega$. We can assume in
this case that $u(0,x)=p_{a}$.

For our computations, $p_{l}$ was obtained by numerically solving
(\ref{eqsta}) thanks to a second order finite elements method with
triangular mesh elements. In the case with Allee effect, $p_{a}$
was computed as the limit of the solution $w(T,x)$, as $T\to
+\infty$ of the initial-value problem
\par\nobreak\noindent  \begin{equation*}
\frac{\partial w}{\partial t}-D \nabla^2 w=w \ g(0,x,w), \ t>0, \
x\in\Omega,
\end{equation*}
with $w(0,x)=2K$ and $\frac{\partial w}{\partial n}(t,x)=0$ over
$\partial \Omega$. The solution $w(T,x)$ was also obtained thanks
to a finite elements method. The convergence of $w(T,x)$ to a
positive equilibrium, as $T\to +\infty$, ensures that the
condition for the existence of $p_{a}$ is fulfilled.

\section*{Acknowledgements} The authors would like to
thank the editor and the anonymous referees for their valuable
suggestions and insightful comments. The numerical computations
were carried out using Comsol Multiphysics$^\circledR$. This study
was supported by the french ``Agence Nationale de la Recherche"
within the project URTICLIM ``Anticipation des effets du
changement climatique sur l'impact \'ecologique et sanitaire
d'insectes forestiers urticants" and by the European Union within
the FP 6 Integrated Project ALARM- Assessing LArge-scale
environmental Risks for biodiversity with tested Methods
(GOCE-CT-2003-506675).


\begin{thebibliography}{99}

\bibitem[Allee(1938)]{allee} Allee WC (1938) The social life of animals. Norton, New York

\bibitem[Amann(1976)]{am2} Amann H (1976) Supersolution, monotone iteration and stability. J Differ Equations 21:367-377

\bibitem[Aronson and Weinberger(1978)]{aw} Aronson DG,  Weinberger HF (1978) Multidimensional nonlinear diffusions arising in population genetics. Adv Math 30:33-76


\bibitem[Berec et al.(2007)]{berec} Berec L, Angulo E, Courchamp F
(2007) Multiple Allee effects and population management. Trends
Ecol Evol 22:185-191

\bibitem[Berestycki and Hamel(2006)]{bh2}
Berestycki H, Hamel F (2006) Fronts and invasions in general
domains. C R Acad Sci Paris Ser I 343:711-716

\bibitem[Berestycki and Lions(1980)]{bl} Berestycki H, Lions P-L (1980) Une m\'ethode locale pour l'\'existence de solutions positives de prol\`emes semi-lin\'eaires elliptiques. J Analyse Math 38:144-187

\bibitem[Berestycki and Rossi(2008)]{br} Berestycki H, Rossi L (2008)
Reaction-diffusion equations for population dynamics with forced
speed I - The case of the whole space. Discret Contin Dyn S, in
press

\bibitem[Berestycki at al.(2005)]{bhr1}  Berestycki H, Hamel F,  Roques L (2005)
  Analysis of the periodically fragmented environment model~: I~- Species persistence. J Math Biol 51:75-113


\bibitem[Berestycki et al.(2007)]{bdnz} Berestycki H, Diekmann O, Nagelkerke CJ, Zegeling PA (2007) Can a species keep pace with a shifting climate?  Bull Math
Biol, submitted


\bibitem[Cantrell and Cosner(2003)]{ccL}  Cantrell RS, Cosner C  (2003) Spatial ecology via reaction-diffusion equations. Series In Mathematical and Computational Biology, John Wiley and Sons, Chichester, Sussex, UK

\bibitem[Chapuisat and Grenier(2005)]{chapuisat}  Chapuisat G,  Grenier E (2005) Existence and nonexistence of traveling wave solutions for a bistable reaction-diffusion equation in an infinite cylinder whose diameter is suddenly increased. Commun Part Differ Eq 30:1805-1816

\bibitem[Deasi and Nelson(2005)]{deasi} Deasi MN, Nelson DR (2005) A quasispecies on a moving oasis. Theor
Popul Biol 67:33-45

\bibitem[Dennis(1989)]{dennis} Dennis B (1989) Allee effects: population growth, critical density, and the chance of extinction. Natural Resource Modeling 3:481-538

\bibitem[Jaeger and Fahrig(2004)]{fahrig2004} Jaeger JAG, Fahrig L (2004) Effects of road fencing on population persistence. Conserv Biol 18:1651-1657

\bibitem[Fife(1979)]{fife} Fife PC (1979) Long-time behavior of solutions of bistable non-linear diffusion equations. Arch Ration Mech An 70:31-46

\bibitem[Fisher(1937)]{fi} Fisher RA (1937)  The wave of advance of advantageous genes. Ann Eugenics~7:355-369


\bibitem[Groom(1998)]{groom} Groom MJ (1998) Allee effects limit population viability of an annual plant. Am Nat 151:487-496

\bibitem[Hilker et al.(2005)]{hilker}
Hilker FM, Lewis MA, Seno H, Langlais M, Malchow H (2005)
Pathogens can slow down or reverse invasion fronts of their hosts.
Biol Invasions 7:817-832

\bibitem[Hurford et al.(2006)]{hurford} Hurford A, Hebblewhite M, Lewis MA (2006) A spatially-explicit
model for the Allee effect: Why wolves recolonize so slowly in
Greater Yellowstone. Theor Popul Biol 70:244-254

\bibitem[IPCC(2007)]{IPCC2007} Intergorvernmental Panel on Climate Change
(2007) Climate change 2007: The physical science basis.
Contribution of working group I to the fourth assessment report of
the intergorvernmental panel on climate change. Summary for
policymakers

\bibitem[Keitt et al.(2001)]{keitt} Keitt TH, Lewis MA, Holt RD (2001) Allee effects, invasion pinning, and species' borders. Am Nat 157:203-216

\bibitem[King and McCabe(2003)]{kmcc}
King JR, McCabe PM (2003) On the Fisher-KPP equation with fast
nonlinear diffusion. Proc R  Soc A Math  Phys Eng Sci
459:2529-2546

\bibitem[Kolmogorov et al.(1937)]{kpp} Kolmogorov AN, Petrovsky IG,  Piskunov NS (1937) Etude de l'\'equation de la diffusion avec croissance de la quantit\'e de mati\`ere et son application \`a un probl\`eme biologique. Bull Univ Etat Moscou, S\'erie Internationale A1:1-26


\bibitem[Lewis and Kareiva(1993)]{lk} Lewis MA, Kareiva P (1993) Allee dynamics and the speed of invading organisms. Theor Popul Biol 43:141-158

\bibitem[Lutscher et al.(2006)]{lutscher} Lutscher F, Lewis MA, McCauley E (2006) Effects of heterogeneity
on spread and persistence in rivers. Bull Math Biol 68:2129-2160


\bibitem[Matano et al.(2006)]{matano} Matano H,  Nakamura K-I,  Lou B (2006)
Periodic traveling waves in a two-dimensional cylinder with
saw-toothed boundary and their homogenization limit. Networks and
Heterogeneous Media 1:537-568

\bibitem[McCarthy(1997)]{mc} McCarthy MA (1997) The Allee effect, finding mates and theoretical models. Ecol Model 103:99-102

\bibitem[Okubo and Levin(2002)]{oku}  Okubo A,  Levin SA (2002) Diffusion and ecological problems - modern perspectives, Second edition, Springer-Verlag

\bibitem[Owen and Lewis(2001)]{owen}
Owen MR, Lewis MA (2001) How predation can slow, stop or reverse a
prey invasion. Bull Math Biol 63:655-684

\bibitem[Pachepsky et al.(2005)]{pach} Pachepsky E, Lutscher F, Nisbet RM, Lewis MA (2005) Persistence,
spread and the drift paradox. Theor Popul Biol 67:61-73

\bibitem[Parmesan(2006)]{pa} Parmesan C (2006) Ecological and evolutionary responses to recent climate change. Annu Rev Ecol Evol Syst 37:637-669

\bibitem[Parmesan and Yohe(2003)]{py} Parmesan C, Yohe G (2003) A globally coherent fingerprint of climate change impacts across natural systems. Nature 421:37-42

\bibitem[Potapov and Lewis(2004)]{pl} Potapov A, Lewis MA (2004) Climate and competition: the effect of moving range boundaries on habitat invisibility. Bull Math Biol 66:975-1008

\bibitem[Robinet et al.(2007a)]{rob} Robinet C, Baier P, Pennerstorfer J, Schopf A, Roques A
(2007a) Modelling the effects of climate change on the potential
feeding activity of \textit{Thaumetopoea pityocampa} (Den. \&
Schiff.) (Lep., Notodontidae) in France. Global Ecol Biogeogr
16:460-471

\bibitem[Robinet et al.(2007b)]{oik} Robinet C, Liebhold A, Gray D (2007b) Variation in developmental time affects mating success and Allee
effects. Oikos 116:1227-1237

\bibitem[Shi and Shivaji(2006)]{shiallee} Shi J,  Shivaji R (2006) Persistence in diffusion models with weak Allee
effect. J Math Biol 52:807-829

\bibitem[Shigesada and Kawasaki(1997)]{sk}  Shigesada N,  Kawasaki, K. (1997) Biological invasions: theory and practice. Oxford Series in Ecology and Evolution, Oxford University Press, Oxford

\bibitem[Stephens and Sutherland(1999)]{steph}  Stephens PA, Sutherland WJ (1999) Consequences
of the Allee effect for behaviour, ecology and conservation.
Trends Ecol Evol 14:401-405

\bibitem[Thomas et al.(2004)]{thomas} Thomas CD, Cameron A, Green RE, Bakkenes M, Beaumont LJ, Collingham YC, Erasmus BFN, de Siqueira MF, GraingerA, HannahL, Hughes L, Huntley B, Jaarsveld AS, Midgley GF, MilesL, Ortega-Huerta MA, Peterson AT, Phillips OL, Williams SE (2004) Extinction risk from climate change. Nature 427:145-148

\bibitem[Tobin et al.(2007)]{tob} Tobin PC, Whitmire SL,  Johnson DM, Bjornstad ON,  Liebhold AM (2007) Invasion speed is affected by geographic variation in the strength of Allee effects. Ecol Lett 10:36-43

\bibitem[Turchin(1998)]{tur} Turchin P (1998) Quantitative analysis of movement: measuring and modeling population redistribution in animals and plants, Sinauer Associates, Sunderland, MA

\bibitem[Veit and Lewis(1996)]{veit} Veit RR, Lewis MA (1996) Dispersal, population growth, and the Allee effect: dynamics of the house finch invasion of eastern North America. Am Nat 148:255-274

\bibitem[Walther et al.(2002)]{walt} Walther GR, Post E, Convey P, Menzel A, Parmesan C, Beebee TJC, Fromentin J-M, Hoegh-Guldberg O,  Bairlein F (2002) Ecological responses to recent climate change. Nature 416:389-395




\end{thebibliography}
\end{document}